\documentclass[10pt,english, a4paper]{article}
\usepackage {amssymb}
\usepackage[intlimits] {amsmath}
\usepackage {babel}
\usepackage {exscale}
\usepackage{epsfig}
\usepackage{color}
\usepackage {theorem}           
\usepackage[latin1]{inputenc}
\usepackage{makeidx}
\usepackage{setspace}
\usepackage{subfigure}
\theoremstyle{break}

\renewcommand{\epsilon}{\varepsilon}

\newcommand{\dx}{\,{\rm{dx}} }

\DeclareMathOperator{\Curl}{Curl}

\newcommand{\norm}[1]{\|#1\|}
\newcommand{\abs}[1]{|#1|}

\newcommand{\R}{\mathbb{R}}

\newcommand{\C}{\mathbb{C}}

\DeclareMathOperator{\GL}{GL}
\DeclareMathOperator{\SO}{SO}

\DeclareMathOperator{\sym}{sym}
\DeclareMathOperator{\Tr}{tr}

\DeclareMathOperator{\Lin}{Lin}

\newcommand{\Mprod}[2]{ {\langle #1 ,#2\rangle} }
\newcommand{\id}{ {1\!\!\!\:1 } }

\DeclareMathOperator{\Det}{det}

\renewcommand{\det}[1]{ {\Det[{#1}]} }
\newcommand{\tr}[1]{ {\Tr \left[{#1}\right]} }

\theoremstyle{break}
\allowdisplaybreaks[1]
\makeatletter
\@addtoreset{equation}{section}

\makeindex 

\def\dist{\hbox{\rm dist}\,}
\def\norm#1{|\!|#1|\!|}
\def\Om{\Omega}

\begin{document}

\title{Counterexamples in the theory of coerciveness for linear elliptic systems related to generalizations of Korn's second inequality}

\author {Patrizio Neff
\thanks {Corresponding author: Patrizio Neff, Lehrstuhl f\"ur Nichtlineare Analysis und Modellierung, Fakult\"at f\"ur Mathematik, Universität Duisburg-Essen, Campus Essen, Thea-Leymann Str. 9, 45127 Essen, Germany, email: patrizio.neff@uni-due.de, Phone +49 201 183 4243, Fax: +49 201 183 4394}\, and
 Waldemar Pompe
 \thanks{Waldemar Pompe, Institute of Mathematics, University of Warsaw, ul. Banacha 2, 02-097 Warszawa, Poland, email: pompe@mimuw.edu.pl. W. Pompe is supported by the Polish Ministry of Science grant no. N N201 397837
(years 2009-2012).} }
\maketitle

\begin{abstract}
 We show that the following generalized
version of Korn's second inequality with nonconstant measurable matrix valued coefficients $P:\Omega\subset\R^3 \to\R^{3\times 3}$
$$\norm{DuP+(DuP)^T}_q+\norm{u}_q\geq c\, \norm{Du}_q\quad
\hbox{for $u\in W_0^{1,q}(\Omega;\R^3), \quad 1<q<\infty$}$$ is in general false, even if $P\in\SO(3)$, while the Legendre-Hadamard condition and ellipticity on $\C^n$ for the 
quadratic form $\abs{Du P+(DuP)^T}^2$ is satisfied.  Thus G\r{a}rding's inequality may be violated for formally positive quadratic forms.
\end{abstract}
\vspace{1cm}


\section{Introduction}
G\r{a}rding's inequality plays a crucial role in the theory of elliptic partial differential equations and systems of equations. In the case of systems related to the linear elasticity second order systems, with which we are concerned, G\r{a}rding's inequality gives sufficient conditions for weak coercivity. More precisely, let $\Omega\subset\R^n$ be a bounded domain and let a mapping $A:\Omega\mapsto \Lin(\R^{n\times n},\R^{n\times n})$ be given. Define a bilinear form
\begin{align}
  a(u,v):=\int\nolimits_\Omega \Mprod{A(x).Du}{Dv}\dx\, ,\quad\quad  u,v\in C_0^\infty(\Omega,\R^n)
  \end{align}
for simplicity without lower order terms. Here, $\Mprod{X}{Y}:=\sum_{i,j=1}^n \, X_{ij} Y_{ij}$ for $X,Y\in \R^{n\times n}$. 

The problem is, under what set of assumptions on $A$ does {\bf weak coercivity} hold, i.e.
\begin{align}
\label{weak_coercivity}
 \exists \, \lambda,\, c>0\;\; \forall \, u\in C_0^\infty(\Omega,\R^n):\quad a(u,u)+\lambda \norm{u}_2^2\ge c\, \norm{Du}_2^2\, .
\end{align}
This is a vector-valued form of G\r{a}rding's inequality \cite{Garding53}. Note that G\r{a}rding's inequality makes a statement about functions with compact support, only. By {\bf strong coercivity} we mean an inequality of the type
\begin{align}
\label{strong_coercivity}
 \exists \, c>0\;\; \forall \, u\in C_0^\infty(\Omega,\R^n):\quad a(u,u) \ge c\,  \norm{Du}_2^2\, .
\end{align}
It is well known \cite[p.74]{Valent88} that G\r{a}rding's inequality \eqref{weak_coercivity} is true, provided that $A$ is uniformly continuous on $\Omega$ (continuous up to the boundary) and $A$ satisfies a uniform Legendre-Hadamard condition
\begin{align}
\label{Legendre_Hadamard}
\exists \, c>0\; \forall\, \xi,\eta\in \R^n\setminus\{0\}:\quad   \Mprod{A(x).(\xi\otimes\eta)}{\xi\otimes\eta}\ge c\, \abs{\xi}^2\, \abs{\eta}^2\, .
\end{align}
It is known that Caccioppoli's inequality, which is an integral inequality estimating the derivatives $Du$ of weak solutions of the corresponding elliptic system in terms of $u$ itself and which is decisive for showing regularity, may break down for merely measurable $A$ satisfying the Legendre-Hadamard ellipticity condition, see \cite{Giaquinta85}, while it is true for uniformly continuous $A$. K. Zhang \cite{Zhang89a,Zhang89b} presented an example such that the mapping $A$ has measurable coefficients, $\Omega\subset\R^3$, satisfying
\begin{align}
    \Mprod{A(x).(\xi\otimes\eta)}{\xi\otimes\eta}=\abs{\xi}^2\, \abs{\eta}^2\,,
\end{align}
but 
\begin{align}
\label{eq:Anonpositive}
\forall\, \lambda>0 \;\; \exists \, u\in C_0^\infty(\Omega,\R^3):\quad \int\nolimits_\Omega \Mprod{A(x).Du}{Du}
+\lambda \, \abs{u}^2\dx<0\, .
\end{align}
However, this negative answer for weak coercivity of the bilinear form $a$ for $A$ with measurable coefficients is using a quadratic form $a$ which is not formally positive. By {\bf formal positivity} we understand 
that there exists some mapping $A$ (which then corresponds to the square-root of $A$ from \eqref{eq:Anonpositive}) such that $a$ can be written as
\begin{align}
 a(u,u)=\int\nolimits_\Omega \, \abs{A(x). Du}^2\dx \ge 0\, .
\end{align} 
For further use let us define for a given continuous mapping $\hat{A}: \R^{n\times n}\mapsto \Lin(\R^{n\times n},\R^{n\times n})$
\begin{align}
   a_{\Phi}(u,v)=\int\nolimits_\Omega \Mprod{\hat{A}(\Phi(x)).Du}{Dv}\dx\, ,\quad u,v\in C_0^\infty(\Omega,\R^n)\,,
\end{align}
where $\Phi\in L^\infty(\Omega,\R^{n\times n})$, at least. For such a structure we have \cite[Th. 6.5.1, p.253]{Morrey66}
\begin{align}
\label{Morrey}
\exists \, \lambda, c>0\quad  \forall\, u\in H_0^1(\Omega):\quad  a_{D\varphi}(u,u)+\lambda(\varphi)\, \norm{u}_2^2 \ge c\, \norm{Du}_2^2
\end{align}
if  $\varphi\in C^1(\overline{\Omega},\R^n)$ and $A(D\varphi(x))$ satisfies the uniform Legendre-Hadamard ellipticity condition. If the bilinear form is formally positive, has variational structure $a(u,u)=\int\nolimits_\Omega \abs{\tilde{A}(x). Du}^2\dx$ and 
satisfies the additional ellipticity condition
\begin{align}
\label{complex_ellipticity_condition}
&\tilde{A}(x)\in \Lin(\C^{n\times n},\C^{n\times n})\, ,\\ \notag
\tilde{A}(x).(\xi\otimes\eta)\neq & 0\quad \text{whenever $\xi,\eta\in\C^n$ and $\xi\neq 0, \eta\neq 0$}\,,
\end{align}
which implies the Legendre-Hadamard ellipticity condition, then one has the stronger inequality for functions without boundary conditions
\begin{align}
\exists \, \lambda, c>0\; \forall\, u\in H^1(\Omega):\quad  a(u,u)+\lambda\, \norm{u}_2^2 \ge c\, \norm{Du}_2^2
\end{align}
if $\tilde{A}$ is uniformly continuous on $\overline{\Omega}$, see \cite{Necas68a,Necas68b}.\\

We are especially interested in such quadratic form that arise in generalizations of Korn's inequality \cite{Ciarlet10,Neff_Pauly_Witsch_cracad11}, namely we consider the formally positive bilinear form
\begin{align}
\label{extended_korn_form}
   a_{P}(u,v):&= \int\nolimits_\Omega \Mprod{\sym (Du\, P)}{\sym(Dv\, P)}\dx=\int\nolimits_\Omega \Mprod{[\sym (Du\, P)]P^T}{Dv}\dx\end{align}
for given $P\in L^\infty(\Omega,\GL^+(3))$. This defines the operator $\hat{A}$ from above via $\hat{A}(P).X:=[\sym(X P)]P^T$. Here, $\sym X:=\frac{1}{2}(X+X^T)$. Bilinear forms having this nonstandard shape are nonetheless ubiquitous, they appear e.g. in micromorphic elasticity models \cite{Klawonn_Neff_Rheinbach_Vanis09,Neff_Forest_jel05,Neff_micromorphic_rse_05}, in geometrically exact formulations of 
plasticity \cite{Neff01c,Neff01d}, in Cosserat models \cite{Neff_Cosserat_plasticity05} or in thin shell 
models \cite{Birsan_Neff_MMS_2013,Birsan_Neff_JElast2013,Neff_plate05_poly}. 

Note that $a_P$ cannot be reduced to a quadratic form of the linearized elastic strains $\sym Du$ for general $P$. Further discussions of coercivity for quadratic forms defined on the linearized strains and lack of coercivity for such models with inhomogeneous material parameters but satisfying the Legendre-Hadmard ellipticity condition can be found in \cite{LeDret87,Zhang2003,Zhang2010,Zhang2011}.

\subsection{The geometrically exact Cosserat model}
In order to see the significance of the new bilinear form \eqref{extended_korn_form} we briefly introduce the variational {\bf geometrically exact Cosserat model}: the goal in this extended continuum model is to find the deformation $\varphi:\Omega\subset\R^3\mapsto\R^3$ and the Cosserat microrotation $\overline{R}:\Omega\subset\R^3\mapsto\SO(3)$
\begin{align}
\label{Cosserat_minimization_energy}
 \int\nolimits_\Omega &W(D\varphi,\overline{R})+W_{\rm curv}(D \overline{R})-\Mprod{f}{\varphi}\dx\mapsto\;\min. \;(\varphi,\overline{R}) \, , \notag\\
   W(D\varphi,\overline{R})&=\mu_e\,\abs{\sym(\overline{R}^TD \varphi -\id)}^2
    +\frac{\lambda_e}{2}\,\tr{\sym(\overline{R}^T D\varphi-\id)}^2\, , \notag\\
    W_{\rm curv}(D\overline{R})&=\mu_e\, \left( \frac{L_c^2}{2}\, \abs{\Curl \overline{R}}^2+ \frac{L_c^q}{q}\abs{\Curl\overline{R}}^q\right)\, .
  \end{align}
 Here, $L_c>0$ is defining an intrinsic length scale in the model, while the elastic Lam\'e coefficients satisfy $\mu_e, 3\lambda_e+2\mu_e>0$.
The model is geometrically exact in the sense that it is invariant under the rigid rotation $(\varphi,\overline{R})\mapsto (Q\varphi,Q\overline{R})$ for any constant $Q\in \SO(3)$. This sets the model apart from linear elasticity. Existence for this model hinges on the coerciveness properties 
of 
\begin{align}
\int\nolimits_\Omega\abs{\overline{R}^T D\varphi+D\varphi^T \overline{R} }^2\dx=\int\nolimits_\Omega\abs{ D\varphi\overline{R}^T+\overline{R}D\varphi^T  }^2\dx=4\, a_{\overline{R}^T}(\varphi,\varphi)
\end{align}
at given rotation tensor $\overline{R}\in\SO(3)$. Our result below shows that coercivity of the Cosserat model w.r.t. deformations $\varphi$ needs some additional smoothness which can be ensured via the curvature contribution $W_{\rm curv}$, see \cite{Neff_curl06}. For applications of the Cosserat model in materials science, we refer to \cite{Neff_Muench_magnetic08,Neff_Muench_transverse_cosserat08,Neff_Muench_simple_shear09}.

For $P=\id$ we obtain
\begin{align}
   a_\id(u,u)=\int\nolimits_\Omega \Mprod{\sym  Du}{\sym Du}\dx\,,
 \end{align}
which is a measure for the linear elastic strain. G\r{a}rding's inequality is then nothing else but a simplified version of 
Korn's second inequality \cite{Gobert62,Nitsche81} on $H_0^1(\Omega)$, i.e.
\begin{align}
  \int\nolimits_\Omega\abs{\sym Du}^2+\abs{u}^2\dx\ge c\, \norm{Du}_2^2\, .
\end{align}
Since the bilinear form $a_P$ satisfies a uniform Legendre-Hadamard ellipticity condition
\begin{align}
    \Mprod{\sym ((\xi\otimes\eta)\, P)}{\sym((\xi\otimes\eta)\, P)}&= \Mprod{\sym (\xi\otimes P^T\eta)}{\sym(\xi\otimes P^T\eta)}\\ \notag
     &\ge \frac{1}{2} \abs{\xi}^2\, \abs{P^T \eta}^2\ge \frac{1}{2}\lambda_{\rm min}(PP^T)\, \abs{\xi}^2\, \abs{\eta}^2
     \end{align}
for $P$ such that $\det{P}\ge \mu >0$ and $P\in C(\overline{\Omega}, \R^{n\times n})$, we infer weak coercivity e.g. from \cite[p.74]{Valent88}, i.e.
\begin{align}
  \exists \, \lambda,c>0\, \forall\, u\in H_0^1(\Omega):\quad  a_{P}(u,u)+\lambda\, \norm{u}_2^2\ge c\,  \norm{Du}_2^2\, .
\end{align}
Due to the special structure of the bilinear form $a_P$ it is easy to see that the ellipticity condition \eqref{complex_ellipticity_condition} is also satisfied and we know furthermore that $a_P$ is strictly coercive, provided that $\det P\ge \mu >0$ and $P\in C(\overline{\Omega}, \R^{n\times n})$, see \cite{Neff_Lankeit_Pauly_cracad13,Neff_Lankeit_Pauly_zamp13,Neff00b,Pompe03}. If $P$ is invertible but merely measurable, then we know that strong coercivity, i.e. Korn's first inequality, is in general not true \cite{Pompe10}. If $P$ is invertible, measurable, symmetric and positive definite, then strict coercivity in $H_0^1(\Omega)$ is obtained, without further smoothness assumptions \cite{Pompe03}. Finally, if $P^{-1}=D\varphi$ for a diffeomorphism $\varphi\in C(\overline{\Omega},\R^n)$ ($D\varphi\in L^\infty$) then strict coercivity is obtained as well by a simple transformation of variables argument.

By and large, $a_P$ from \eqref{extended_korn_form} is not strictly coercive if $P$ is only invertible and measurable. Nevertheless, weak coercivity for $a_P$, i.e., the generalization of Korn's second inequality,  could still be true. 

However, in this contribution we show by way of counterexamples that weak coercivity
\begin{align}
   \exists \, \lambda,c>0\, \forall\, u\in H_0^1\Omega):\quad  a_{P}(u,u)+\lambda\, \norm{u}_2^2\ge c\, \norm{Du}_2^2\, 
\end{align}
fails in general for $P$ invertible and measurable. We generalize our counterexamples in the obvious way to the $L^q$-setting, i.e.
\begin{align}
\label{korn_weak_coercivity_q}
  \exists \, \lambda,c>0\,\; \forall\, u\in W_0^{1,q}(\Omega,\R^n):\quad  \int\nolimits_\Omega \abs{\sym Du P}^q+\lambda\, \abs{u}^q \dx \ge c\,  \norm{Du}_q^q\, 
\end{align}
does not hold for $q>1$.

\section{Main Part}
For simplicity, we restrict ourselves to three space dimensions. From now on let $\Omega$ be an open, bounded set in $\R^3$ and let 
$P\colon\,\Omega\subset\R^3\to \R^{3\times 3}$ with $P\in L^\infty(\Omega)$ and
$\det{P(x)}\geq \mu>0$ be given. Assume moreover $q>1$. 

Inequality \eqref{korn_weak_coercivity_q} is equivalent to  
\begin{align}
\label{korn_weak_coercivity_pompe}
 \exists\; c>0\;\; \forall\;u\in W_0^{1,q}(\Omega;\R^3):\quad     \norm{DuP+(DuP)^T}_q+\norm{u}_q\geq c\norm{Du}_q\, .
 \end{align}

It is known that inequality \eqref{korn_weak_coercivity_pompe} holds, if we additionally assume that 
$P\in C(\overline\Om)$ \cite{Necas68a,Necas68b}. This is a generalization of Korn's
second inequality. However, 
we shall show that this inequality is not in general true with noncontinuous, 
bounded coefficients $P$, under the structural assumption 
$P\in SO(3)$. We even show that the following weaker 
inequality is not valid in this case:
\begin{align}
\label{pompe2}
\norm{DuP+(DuP)^T}_q+\norm{u}_\infty\geq c\norm{Du}_q\quad
\hbox{for $u\in W_0^{1,q}(\Omega;\R^3)\cap L^\infty(\Om;\R^3)$}\,.
\end{align}

We present two counterexamples to inequality \eqref{pompe2}
The first one assumes that the coefficients $P(x)$ are bounded and satisfy
$\det{P(x)}=1$ a.e. on $\Om$, while the second counterexample 
assumes more about the
coefficients: $P\in SO(3)$. The construction in the second case is based on the
following result by A. Cellina and S. Perrotta \cite{Cellina95}: 
{\it If $\Om$ is an open, bounded set in $\R^3$, then there exists a mapping
$u\in W_0^{1,\infty}(\Om;\R^3)$ such that $Du(x)\in O(3)$ a.e. on $\Om$.}

Even if the
counterexample in the second case provides also a 
counterexample in the first one,
we present a more elementary construction in the 
first case, which does not require the strong result 
of A. Cellina and S. Perrotta \cite{Cellina95}.
The ideas of our constructions are similar to the constructions presented by the second author in \cite{Pompe03,Pompe10}.

Our method of construction is direct and yields $P$ having a finite number of elements. This result - due to a complicated structure of the rank-one connections - would be hard to obtain by the convex integration method.

To deal better with constants, we use the following 
definition of the $L^q$-norm of a
mapping $P\colon\, \Om\to\R^{3\times 3}$. For 
$$P(x)=\begin{pmatrix}
                   p_1(x)\\
                   p_2(x)\\
                   p_3(x)
                   \end{pmatrix} 
\quad (x\in\Om)$$
define
$$\norm{P}_q^q=\int\nolimits_\Om(|p_1(x)|^q+|p_2(x)|^q+|p_3(x)|^q)\,dx\,,$$
where $|p_i|$ denotes the Euclidean norm of the vector $p_i\in\R^3$.

{\bf Theorem 1.}

For each $q>1$ and any open, bounded subset $\Omega$ of $\R^3$,
there exist $P\in L^\infty(\Om;\R^{3\times 3})$ with
$\det P(x)=1$ and a sequence 
$u_n\in W_0^{1,q}(\Om;\R^3)\cap L^\infty(\Om;\R^3)$, such that:

\smallskip
(a) $Du_nP+(Du_nP)^T=0$ on the set $\Omega$,

\smallskip
(b) $\norm{Du_n}_q= 2^{1/q}$,

\smallskip
(c) $\norm{u_n}_\infty\to0$ as $n\to\infty$.

\medskip\goodbreak
{\bf Proof}

Let $\Om_1,\Om_2,\ldots$ be open, disjoint subsets of the set $\Om$, such that
the set $\Om\setminus(\Om_1\cup\Om_2\cup\ldots)$ has the measure zero.
Let moreover $R$ be a fixed rotation in $\R^3$ such that $Re_i\neq \pm e_j$
for all $i,j\in\{1,2,3\}$. 
On each of the sets $\Om_n$ we construct two Vitali coverings: one of them
with cubes $Q_{ni}$ $(i=1,2,\ldots)$, whose edges are parallel to the vectors $e_1$,
$e_2$ and $e_3$ and the other with cubes
$S_{nj}$ $(j=1,2,\ldots)$, whose edges are parallel to the
vectors $Re_1$, $Re_2$ and $Re_3$.

Therefore we have $\Om_n=Q_{n1}\cup Q_{n2}\cup\ldots$, where the interiors of the cubes
$Q_{ni}$ $(i=1,2,\ldots)$ are disjoint, and similarly $\Om_n=S_{n1}\cup S_{n2}\cup\ldots$, 
where the interiors of the cubes
$S_{nj}$  $(j=1,2,\ldots)$ are disjoint.
We may moveover assume that the length of the edge of each cube
$Q_{ni}$ and $S_{nj}$ is at most $\displaystyle{2|\Om_n|^{1/q}\over n}$.

Now on each of the sets $\Om_n$ define two mappings: $u_n^1,\;u_n^2\colon\,\Om_n\to\R$ 
as follows:
$$u_n^1(x)=\dist(x,\partial Q_{ni})\quad\hbox{for $x\in Q_{ni}$}\quad(i=1,2,\ldots)$$
and
$$u_n^2(x)=\dist(x,\partial S_{nj})\quad\hbox{for $x\in S_{nj}$}\quad(j=1,2,\ldots)\,.$$
Then $u_n^1,\;u_n^2\in W_0^{1,\infty}(\Om_n)$ and $|Du_n^1(x)|=|Du_n^2(x)|=1$ a.e. on $\Om_n$.
Moreover, we have
$$|u_n^1(x)|\leq {|\Om_n|^{1/q}\over n}\quad\hbox{and}\quad 
|u_n^2(x)|\leq {|\Om_n|^{1/q}\over n}
\quad(x\in\Om_n)\,.\leqno(*)$$

Since $Re_i\neq \pm e_j$, the vectors $Du_n^1(x)$, 
$Du_n^2(x)$ are not parallel and therefore
they span a 2-dimensional plane $\pi(x)$ in $\R^3$. Let $v_n(x)$ 
be the vector orthogonal to this plane, such that
$${\rm{det}}\begin{pmatrix} 
               -Du_n^2(x)\\ Du_n^1(x) \\ v_n(x)
\end{pmatrix}=1
\quad(x\in\Om_n)\,.$$
Then since $Re_i\neq \pm e_j$, we obtain $c_1<|v_n^3(x)|<c_2$, where the constants
$c_1$ and $c_2$ are positive and depend only on $R$. 

Define
$$P_n(x)=\begin{pmatrix}
                      -Du_n^2(x)\\ 
                       Du_n^1(x)\\ 
                       v_n(x)
                       \end{pmatrix}^{-1}\quad\hbox{for $x\in\Om_n$}$$
and $P_n(x)=0$ for $x\in\Om\setminus\Om_n$. 
Then $P_n\in L^\infty(\Om_n;\R^{3\times 3})$. Finally let
$$P(x)=\displaystyle\sum_{n=1}^\infty P_n(x)\,.$$ 
Since the supports of the mappings
$P_n(x)$ are disjoint, the above sum is actually a single summand for almost 
each $x\in\Om$ and therefore $P(x)=P_n(x)$
for a.e. $x\in\Om_n$. It follows therefore that 
$P\in L^\infty(\Om;\R^{3\times 3})$.
It is also clear that $\det{P(x)}=1$.

Now define the mappings $u_n\colon\,\Om\to\R^3$ with
$$u_n(x)={1\over |\Om_n|^{1/q}}\,(u_n^1(x),u_n^2(x),0)\quad
\hbox{for $x\in\Om_n$}\,,$$
and $u_n(x)=0$ for $x\in\Om\setminus\Om_n$. Then using $(*)$, we obtain
$$|u_n(x)|\leq {1\over |\Om_n|^{1/q}}\cdot(|u_n^1(x)|+|u_n^2(x)|)\leq 
{2\over n}\quad(x\in\Om)\,.$$
Hence $u_n\in  W_0^{1,\infty}(\Om)$ and the property (c) holds. Moreover,
for $x\in\Om_n$ we have
$$Du_n(x)P(x)=Du_n(x)P_n(x)=
\begin{pmatrix}
0 & |\Om_n|^{-1/q} & 0 \\
-|\Om_n|^{-1/q} & 0 & 0 \\
0 & 0 & 0 
\end{pmatrix}\,,$$
and $Du_n(x)P(x)=0$ for $x\in\Om\setminus\Om_n$.
This shows property (a). Finally, to see property (b) note that
$$\norm{Du_n}_q^q={1\over |\Om_n|}\int\nolimits_{\Om_n} 
\left(|Du_n^1(x)|^q+|Du_n^2(x)|^q\right)\,dx=2\,,$$
and the conclusion (b) follows.

\goodbreak
{\bf Theorem 2.}

For each $q>1$ and any open, bounded subset $\Omega$ of $\R^3$,
there exist a mapping $P\colon\,\Om\to\R^{3\times 3}$ with
$P(x)\in SO(3)$ and a sequence 
$u_n\in W_0^{1,q}(\Om;\R^3)\cap L^\infty(\Om,\R^3)$, such that:

\smallskip
(a) $Du_nP+(Du_nP)^T=0$ for $x\in\Omega$,

\smallskip
(b) $\norm{Du_n}_q= 2^{1/q}$,

\smallskip
(c) $\norm{u_n}_\infty\to0$ as $n\to\infty$.

\medskip\goodbreak
{\bf Proof}

The construction is based on the following 
result of A. Cellina and S. Perrotta \cite{Cellina95}:
{\it If $\Om$ is an open, bounded set in $\R^3$, then there exists a mapping
$u\in W_0^{1,\infty}(\Om;\R^3)$ such that $Du(x)\in O(3)$ a.e. on $\Om$.}

Let $\Om$ be represented, up to a set of measure 0, 
by a union of disjoint open sets $\Om_1,\Om_2,\ldots$.
On each of the sets $\Om_n$ we construct a Vitali covering
with cubes $Q_{ni}$ $(i=1,2,\ldots)$.
Therefore we have $\Om_n=Q_{n1}\cup Q_{n2}\cup\ldots$, 
where the interiors of the cubes
$Q_{ni}$ $(i=1,2,\ldots)$ are disjoint.
We may moveover assume that the length of the edge of each cube
$Q_{ni}$ $(i=1,2,\ldots)$ is at most $\frac{1}{n}|\Om_n|^{1/q}$.

Let $Q$ be a unit cube and let 
$v\in W_0^{1,\infty}(Q;\R^3)$ such that $Dv(x)\in O(3)$ a.e. on $Q$.
Let $u(x)=(v_1(x),v_2(x),0)$. Then $Du_1(x)$, $Du_2(x)$
have the length 1 and are orthogonal for a.e. $x\in Q$.

Set $c=\norm{u}_\infty$. Then scaling and translating the mapping $u$, we 
construct on each cube $Q_{ni}$ a mapping 
$u_{ni}\in W_0^{1,\infty}(Q_{ni};\R^3)$
with 
$$Du_{ni}(x)=\begin{pmatrix}
                             Du_{ni}^1(x)\\
                             Du_{ni}^2(x) \\ 
                             0
                             \end{pmatrix}
\quad \hbox{for a.e. $x\in Q_{ni}$}\,,$$
where $Du_{ni}^1(x)$ and $Du_{ni}^2(x)$ have length $1$ and are orthogonal for
a.e. $x\in Q_{ni}$ and
$$|u_{ni}(x)|\leq c\cdot {|\Om_n|^{1/q}\over n}\quad (x\in Q_{ni})\,.$$

Extend the mappings $u_{ni}$ $(i=1,2,\ldots)$ from $Q_{ni}$ to $\Om$ 
by setting $u_{ni}(x)=0$ on $\Om\setminus Q_{ni}$ and define
$$u_n(x)={1\over |\Om_n|^{1/q}}\sum_{i=1}^\infty u_{ni}(x)\quad(x\in\Om)\,.$$
Then for a.e. $x\in Q_{ni}$ we have 
$$u_n(x)={1\over|\Om_n|^{1/q}}\, u_{ni}(x)\,.$$
Thus
$$|u_n(x)|={1\over|\Om_n|^{1/q}}\, |u_{ni}(x)|\leq 
{c\over n}\quad (x\in\Om)\,.$$
This shows the property (c).

Now define $P(x)$ as follows. If $x\in Q_{ni}$, then define $P(x)$ such that
$P(x)^T$ is the rotation, which takes the vectors
$Du_{ni}^1(x)$, $Du_{ni}^2(x)$ to the vectors 
$(0,1,0)$, $(-1,0,0)$, respectively
(such a rotation exists, since the vectors 
$Du_{ni}^1(x)$, $Du_{ni}^2(x)$ have the length $1$
and are orthogonal). Then we have
$$Du_n(x)P(x)=
\begin{pmatrix}
0 & |\Om_n|^{-1/q} & 0 \\
-|\Om_n|^{-1/q} & 0 & 0 \\
0 & 0 & 0 
\end{pmatrix} \quad\hbox{for $x\in\Om_n$}\,,$$
and $Du_n(x)P(x)=0$ for $x\in\Om\setminus\Om_n$.
This shows property (a). Finally, to see property (b) note that
$$\norm{Du_n}_q^q=\int\nolimits_{\Om_n} \left(|Du_n^1(x)|^q+|Du_n^2(x)|^q\right)\,dx=
\int\nolimits_{\Om_n} {2\over |\Om_n|}\,dx=2\,,$$
and the conclusion (b) follows.

\bibliographystyle{plain} 
{\footnotesize
\bibliography{literatur1}
}

\end{document}